\documentclass[twosided,reqno]{amsart}

\usepackage{amsmath}
\usepackage{amssymb}
\usepackage{amsthm}

\newcommand{\iz}{{\int_{{\mathbb{Z}}_p}}}

\theoremstyle{theorem}
\newtheorem{theorem}{\scshape Theorem }[section]

\newtheorem{corollary}[theorem]{\scshape Corollary}

\theoremstyle{definition}

\numberwithin{equation}{section}

\begin{document}

\title[On degenerate Carlitz $q$-Bernoulli polynomials]{On degenerate Carlitz $q$-Bernoulli polynomials}

\author{Taekyun Kim}
\address{Department of Mathematics, Kwangwoon University, Seoul 139-701, Republic of Korea.}
\email{tkkim@kw.ac.kr}


\maketitle

\begin{abstract}
In this paper, we consider the degenerate Carlitz $q$-Bernoulli numbers and polynomials and we investigate some properties of those polynomials.
\end{abstract}

\section{Introduction}

Throughout this paper, ${\mathbb{Z}}_p$, ${\mathbb{Q}}_p$ and ${\mathbb{C}}_p$ will, respectively, denote the ring of $p$-adic rational integers, the field of $p$-adic rational numbers and the completion of the algebraic closure of ${\mathbb{Q}}_p$. Let $\nu_p$ be the normalized exponential valuation of ${\mathbb{C}}_p$ with ${p}_p=p^{-\nu_p(p)}=\frac{1}{p}$. When one talks of $q$-extension, $q$ is variously considered as an indeterminate, a complex number $q\in{\mathbb{C}}$, or $p$-adic number $q\in{\mathbb{C}}_p$. If $q\in{\mathbb{C}}$, we assume that $|q|<1$. If $q\in{\mathbb{C}}_p$, we assume $|q-1|_p<p^{-\frac{1}{p-1}}$ so that $q^x=\exp(x\log q)$ for $|x|_p<1$. We use the notation $[x]_q=\frac{1-q^x}{1-q}$. Note that $\lim_{q\rightarrow 1}=x$.

In \cite{03},  L. Carlitz considered $q$-Bernoulli numbers as follows:
\begin{equation}\label{1}
\beta_{0,q}=1,~q(q\beta_q+1)^n-\beta_{n,q}=
\begin{cases}
1,&if~n=1,\\
0, &if~n>1,
\end{cases}
\end{equation}
withe the usual convention about replacing $\beta_1 ^n$ by $\beta_{n,q}$. The $q$-Bernoulli polynomials are defined by
\begin{equation}\label{2}
\beta_{n,q} (x)=\sum_{l=0} ^n \binom{n}{l}\beta_{l,q} q^l[x]_q ^{n-l},{\text{ (see \cite{06})}}.
\end{equation}

In \cite{02}, L. Carlitz defined the degenerate Bernoulli polynomials which are given by the generating function to be
\begin{equation}\label{3}
\frac{t}{(1+\lambda t)^{\frac{1}{\lambda}}-1}(1+\lambda t)^{\frac{x}{\lambda}}=\sum_{n=0} ^{\infty}\beta_n(x|\lambda)\frac{t^n}{n!},{\text{ (see \cite{02,03})}}.
\end{equation}
When $x=0$, $\beta_n(x)(\lambda)=\beta_n(0|\lambda)$ are called the {\it{degenerate Bernoulli numbers}}. Note that $\lim_{\lambda\rightarrow 0}\beta_n(x|\lambda)=B_n(x)$, where $B_n(x)$ are ordinary Bernoulli polynomials.

Let $UD({\mathbb{Z}}_p)$ be the space of uniformly differentiable function on ${\mathbb{Z}}_p$. For $f\in UD({\mathbb{Z}}_p)$, the {\it{$p$-adic $q$-integral on ${\mathbb{Z}}_p$}} is defined by Kim to be
\begin{equation}\label{4}
I_q(f)=\iz f(x)d\mu_q(x)=\lim_{N\rightarrow \infty}\frac{1}{\left[p^N\right]_q}\sum_{x=0} ^{p^N-1}f(x)q^x,{\text{ (see \cite{06})}}.
\end{equation}
The Carlitz's $q$-Bernoulli polynomials can be represented by $p$-adic $q$-integral on ${\mathbb{Z}}_p$ as follows:
\begin{equation}\label{5}
\iz[x+y]_q ^nd\mu_q(y)=\beta_{n,q} (x),~(n\geq 0).
\end{equation}
Thus, by \eqref{4}, we get
\begin{equation}\label{6}
\iz e^{[x+y]_q t}d\mu_q(y)=\sum_{n=0} ^{\infty}\beta_{n,q}(x)\frac{t^n}{n!},{\text{ (see [1-8])}}.
\end{equation}
From \eqref{6}, we can derive the following equation:
\begin{equation}\label{7}
\beta_{n,q}(x)=\frac{1}{(1-q)^m}\sum_{j=0} ^m \binom{m}{j}(-1)^jq^{ix}\frac{j+1}{[j+1]_q},~(m\geq0).
\end{equation}
In this paper, we consider the degenerate Carlitz $q$-Bernoulli polynomials and numbers and we investigate some properties of those polynomials.

\section{Degenerate Carlitz $q$-Bernoulli numbers and polynomials}

In this section, we assume that $\lambda,t\in{\mathbb{C}}_p$ with $|\lambda t|_p<p^{-\frac{1}{p-1}}$. In the viewpoint of \eqref{3}, we consider the {\it{degenerate Carlitz $q$-Bernoulli polynomials}} which are given by the generating function to be
\begin{equation}\label{8}
\iz(1+\lambda t)^{\frac{1}{\lambda}[x+y]_q}d\mu_q(y)=\sum_{n=0} ^{\infty} \beta_{n,q}(x|\lambda)\frac{t^n}{n!}.
\end{equation}
When $x=0$, $\beta_{n,q} (\lambda)=\beta_{n,q} (0|\lambda)$ are called the {\it{degenerate Carlitz $q$-Bernoulli numbers}}.

Now, we observe that
\begin{equation}\label{9}
\begin{split}
\iz(1+\lambda t)^{\frac{[x+y]_q}{\lambda}}d\mu_q(y)=&\sum_{n=0} ^{\infty} \iz\binom{\frac{[x+y]_q}{\lambda}}{n}d\mu_q(y)\lambda^nt^n\\
=&\sum_{n=0} ^{\infty}\iz\left(\frac{[x+y]_q}{\lambda}\right)_nd\mu_q(y)\lambda^n\frac{t^n}{n!},
\end{split}
\end{equation}
where $\left(\frac{[x+y]_q}{\lambda}\right)_n=\frac{[x+y]_q}{\lambda} \times\left(\frac{[x+y]_q}{\lambda}-1\right)\times\cdots\times\left(\frac{[x+y]_q}{\lambda}-n+1\right)$.

Now, we define $[x+y]_{n,\lambda}$ as
\begin{equation}\label{10}
[x+y]_{n,\lambda}=[x+y]_q\left([x+y]_q-\lambda\right)\cdots\left([x+y]_q-(n-1)\lambda\right),~(n\geq0).
\end{equation}
Therefore, by \eqref{8}, \eqref{9} and \eqref{10}, we obtain the following theorem.
\begin{theorem}\label{thm1}
For $n \geq 0$, we have
\begin{equation*}
\iz[x+y]_{n,\lambda}d\mu_q(y)=\beta_{n,q}(x|\lambda).
\end{equation*}
\end{theorem}

Let $S_1(n,m)$ be the {\it{Stirling number of the first kind}} which is defined by $(x)_n=\sum_{l=0} ^n S_1(n,l)x^l$, $(n\geq0)$. Then, by \eqref{9}, we get
\begin{equation}\label{11}
\begin{split}
\iz\left(\frac{[x+y]_q}{\lambda}\right)_nd\mu_q(y)=&\sum_{l=0} ^n S_1(n,l)\lambda^{-l}\iz[x+y]_q ^ld\mu_q(y) \\
=&\sum_{l=0} ^n S_1(n,l)\lambda^{-l}\beta_{l,} (x).
\end{split}
\end{equation}
Therefore, by \eqref{9} and \eqref{11}, we obtain the following theorem.
\begin{theorem}\label{thm2}
For $n\geq 0$, we have
\begin{equation*}
\beta_{n,q}(x|\lambda)=\sum_{l=0} ^n S_1(n,l)\lambda^{n-l}\beta_{l,q}(x).
\end{equation*}
\end{theorem}
Note that $\lim_{\lambda\rightarrow 0}\beta_{n,q}(x|\lambda)=\beta_{n,q}(x)$.
\begin{corollary}\label{coro3}
For $n \geq 0$, we have
\begin{equation*}
\beta_{n,q}(x|\lambda)=\sum_{l=0} ^n\sum_{j=0} ^l\frac{S_1(n,l)}{(1-q)^l}\binom{l}{j}(-1)^jq^{jx}\frac{j+1}{[j+1]_q}\lambda^{n-l}.
\end{equation*}
\end{corollary}
We observe that
\begin{equation}\label{12}
\begin{split}
(1+\lambda t)^{\frac{[x+y]_q}{\lambda}}=&e^{\frac{[x+y]_q}{\lambda}\log(1+\lambda t)}=\sum_{n=0} ^{\infty}\left(\frac{[x+y]_q}{\lambda}\right)^n\frac{1}{n!}\left(\log(1+\lambda t)\right)^n\\
=&\sum_{m=0} ^{\infty}\left(\frac{[x+y]_q}{\lambda}\right)^m\frac{1}{m!}m!\sum_{n=m} ^{\infty}S_1(n,m)\frac{\lambda^nt^n}{n!}\\
=&\sum_{n=0} ^{\infty}\left(\sum_{m=0} ^n \lambda^{n-m}S_1(n,m)[x+y]_q ^m\right)\frac{t^n}{n!}.
\end{split}
\end{equation}
Thus, by \eqref{12}, we get
\begin{equation*}
\begin{split}
\iz(1+\lambda t)^{\frac{[x+y]_q}{\lambda}}d\mu_q(y)=&\sum_{n=0} ^{\infty}\left(\sum_{m=0} ^n \lambda^{n-m}S_1(n,m)\iz[x+y]_q ^md\mu_q(x)\right)\frac{t^n}{n!}\\
=&\sum_{n=0} ^{\infty}\left(\sum_{m=0} ^n \lambda^{n-m}S_1(n,m)\beta_{m,q}(x)\right)\frac{t^n}{n!}.
\end{split}
\end{equation*}
Replacing $t$ by $\frac{1}{\lambda}\left(e^{\lambda t}-1\right)$ in \eqref{8}, we get
\begin{equation}\label{13}
\begin{split}
\iz e^{[x+y]_qt}d\mu_q(y)=&\sum_{m=0} ^{\infty}\beta_{m,q}(x|\lambda)\frac{1}{m!}\frac{1}{\lambda^m}\left(e^{\lambda t}-1\right)^m\\
=&\sum_{m=0} ^{\infty}\beta_{m,q}(x|\lambda)\lambda^{-m}\sum_{n=m} ^{\infty}S_2(n,m)\frac{\lambda^nt^n}{n!}\\
=&\sum_{n=0} ^{\infty}\left(\sum_{m=0} ^n \beta_{m,q}(x|\lambda)\lambda^{n-m}S_2(n,m)\right)\frac{t^n}{n!},
\end{split}
\end{equation}
where $S_2(n,m)$ is the Stirling numbers of the second kind.

We note that the left hand side of \eqref{13} is given by
\begin{equation}\label{14}
\begin{split}
\iz e^{[x+y]_qt}d\mu_q(y)=&\sum_{n=0} ^{\infty}\iz [x+y]_q ^n d\mu_q(y)\frac{t^n}{n!}\\
=&\sum_{n=0} ^{\infty}\beta_{n,q}(x)\frac{t^n}{n!}.
\end{split}
\end{equation}
Therefore, by \eqref{13} and \eqref{14}, we obtain the following theorem.
\begin{theorem}\label{thm4}
For $n \geq 0$, we have
\begin{equation*}
\beta_{n,q}(x)=\sum_{m=0} ^n \beta_{m,q} (x|\lambda)\lambda^{n-m}S_2(n,m).
\end{equation*}
\end{theorem}

Note that
\begin{equation}\label{15}
\begin{split}
(1+\lambda t)^{\frac{[x+y]_q}{\lambda}}=&(1+\lambda t)^{\frac{[x]_q}{\lambda}}(1+\lambda t)^{\frac{q^x[y]_q}{\lambda}}\\
=&\left(\sum_{m=0} ^{\infty}[x]_{m,\lambda}\frac{t^m}{m!}\right)\left(\sum_{l=0} ^{\infty}\frac{q^{lx}}{\lambda^l}\frac{[y]_q ^l \left(\log(1+\lambda t)\right)^l}{l!}\right)\\
=&\left(\sum_{m=0} ^{\infty}[x]_{m,\lambda}\frac{t^m}{m!}\right)\left(\sum_{k=0} ^{\infty}\left(\sum_{l=0} ^k \lambda^{k-l}q^{lx}[y]_q ^lS_1(k,l)\right)\frac{t^k}{k!}\right)\\
=&\sum_{n=0} ^{\infty}\left(\sum_{k=0} ^n \sum_{l=0} ^k[x]_{n-k,\lambda}\lambda^{k-l}q^{lx}[y]_q ^lS_1(k,l)\binom{n}{k}\right)\frac{t^n}{n!}.
\end{split}
\end{equation}
Thus, by \eqref{15}, we get
\begin{equation}\label{16}
\begin{split}
\sum_{n=0} ^{\infty}\beta_{n,q}(x|\lambda)\frac{t^n}{n!}=&\sum_{n=0} ^{\infty}\left(\sum_{k=0} ^n \sum_{l=0} ^k[x]_{n-k,\lambda}\lambda^{k-l}q^{lx}\iz [y]_q ^ld\mu_q(y)S_1(k,l)\binom{n}{k}\right)\frac{t^n}{n!}\\
=&\sum_{n=0} ^{\infty}\left(\sum_{k=0} ^n \sum_{l=0} ^k\binom{n}{k}[x]_{n-k,\lambda}\lambda^{k-l}q^{lx}\beta_{l,q}S_1(k,l)\right)\frac{t^n}{n!}.
\end{split}
\end{equation}
Therefore, by \eqref{10}, we obtain the following theorem.
\begin{theorem}\label{thm5}
For $n \geq 0$, we have
\begin{equation*}
\beta_{n,q}(x|\lambda)=\sum_{k=0} ^n \sum_{l=0} ^k\binom{n}{k}[x]_{n-k,\lambda}\lambda^{k-l}q^{lx}S_1(k,l)\beta_{l,q}.
\end{equation*}
\end{theorem}

For $r\in{\mathbb{N}}$, we define the {\it{degenerate Carlitz $q$-Bernoulli polynomials of order $r$}} as follows:
\begin{equation}\label{17}
\begin{split}
&\iz\cdots\iz (1+\lambda t)^{\frac{[x_1+\cdots+x_r+x]_q}{\lambda}}d\mu_1(x_1)\cdots d\mu_q(x_r)\\
=&\sum_{n=0} ^{\infty}\beta_{n,q} ^{(r)}(x|\lambda)\frac{t^n}{n!}.
\end{split}
\end{equation}
We observe that
\begin{equation}\label{18}
\begin{split}
&\iz\cdots\iz(1+\lambda t)^{\frac{[x_1+\cdots+x_r+x]}{\lambda}}d\mu_q(x_1)\cdots d\mu_q(x_r)\\
=&\sum_{m=0} ^{\infty}\lambda^{-m}\iz\cdots\iz[x_1+\cdots+x_r+x]_q ^m d\mu_q(x_1)\cdots d\mu_q(x_r)\frac{1}{m!}\left(\log(1+\lambda t)\right)^m\\
=&\sum_{m=0} ^{\infty}\beta_{m,q} ^{(r)}(x)\lambda^{-m}\sum_{n=m} ^{\infty}S_1(n,m)\frac{\lambda^n}{n!}t^n\\
=&\sum_{n=0} ^{\infty}\left(\sum_{m=0} ^n \lambda^{n-m}\beta_{m,q} ^{(r)}S_1(n,m)\right)\frac{t^n}{n!},
\end{split}
\end{equation}
where $\beta_{m,q} ^{(r)}$ are the Carlitz $q$-Bernoulli numbers of order $r$.

Therefore, by \eqref{17} and \eqref{18}, we obtain the following theorem.
\begin{theorem}\label{thm6}
For $n \geq 0$, we have
\begin{equation*}
\beta_{n,q} ^{(r)}=\sum_{m=0} ^n \lambda^{n-m}\beta_{m,q} ^{(r)}S_1(n,m).
\end{equation*}
\end{theorem}

Replacing $t$ by $\frac{1}{\lambda}\left(e^{\lambda t}-1\right)$ in \eqref{17}, we have
\begin{equation}\label{19}
\begin{split}
&\iz\cdots\iz e^{[x_1+\cdots+x_r+x]_q t}d\mu_q(x_1)\cdots d\mu_q(x_r)\\
=&\sum_{m=0} ^{\infty}\beta_{m,q} ^{(r)}(x|\lambda)\frac{1}{m!}\lambda^{-m}\left(e^{\lambda t}-1\right)^m\\
=&\sum_{m=0} ^{\infty}\beta_{m,q} ^{(r)}(x|\lambda)\lambda^{-m}\sum_{n=m} ^{\infty}S_2(n,m)\frac{\lambda^nt^n}{n!}\\
=&\sum_{n=0} ^{\infty}\left(\sum_{m=0} ^n \lambda^{n-m}\beta_{m,q} ^{(r)}(x|\lambda)S_2(n,m)\right)\frac{t^n}{n!}.
\end{split}
\end{equation}
The left hand side of \eqref{19} is given by
\begin{equation}\label{20}
\begin{split}
&\iz\cdots\iz e^{[x_1+\cdots+x_r+x]_qt}d\mu_q(x_1)\cdots d\mu_q(x_r)\\
=&\sum_{n=0} ^{\infty}\beta_{n,q} ^{(r)}(x)\frac{t^n}{n!}.
\end{split}
\end{equation}
By comparing the coefficients on the both sides of \eqref{19} and \eqref{20}, we obtain the following theorem.
\begin{theorem}\label{thm7}
For $n\geq0$, we have
\begin{equation*}
\beta_{m,q} ^{(r)}(x)=\sum_{m=0} ^n \lambda^{n-m}S_2(n,m)\beta_{m,q} ^{(r)}(x|\lambda).
\end{equation*}
\end{theorem}
We recall that
\begin{equation*}
\begin{split}
\iz f(x)d\mu_q(x)=&\lim_{N\rightarrow \infty}\frac{1}{[p^N]_q}\sum_{x=0} ^{p^N-1}f(x)q^x\\
=&\lim_{N\rightarrow \infty}\frac{1}{[dp^N]_q}\sum_{x=0} ^{dp^N-1}f(x)q^x,
\end{split}
\end{equation*}
where $d\in{\mathbb{N}}$ and $f\in UD({\mathbb{Z}_p})$.

Now, we observe that
\begin{equation}\label{21}
\beta_{n,q}(x|\lambda)=\sum_{l=0} ^n S_1(n,l)\lambda^{n-l}\iz[x+y]_q ^l d\mu_q(y),
\end{equation}
and
\begin{equation}\label{22}
\begin{split}
\iz[x+y]_q ^ld\mu_q(y)=&\frac{1}{[m]_q}\sum_{i=0} ^{m-1}q^i[m]_q ^l\iz\left[\frac{x+i}{m}+y\right]_{q^m} ^ld\mu_{q^m}(y)\\
=&[m]_q ^{l-1}\sum_{i=0} ^{m-1}q^i\beta_{l,q^m}\left(\frac{x+i}{m}\right),
\end{split}
\end{equation}
where $l \in {\mathbb{Z}}_{\geq 0}$ and $m\in{\mathbb{N}}$.

Therefore, by \eqref{21} and \eqref{22}, we obtain the following theorem.
\begin{theorem}\label{thm8}
For $n \geq {\mathbb{N}}_{\geq0}$, $m\in{\mathbb{N}}$, we have
\begin{equation*}
\beta_{n,q}(x|\lambda)=\sum_{l=0} ^n \sum_{i=0} ^{m-1} S_1(n,l)\lambda^{n-l}[m]_q ^{l-1}q^i\beta_{l,q^m}\left(\frac{x+i}{m}\right).
\end{equation*}
\end{theorem}

From \eqref{4}, we note that
\begin{equation}\label{23}
qI_q(f_1)-I_q(f)=(q-1)f(0)+\frac{q -1}{\log q}f^{'}(0),
\end{equation}
where $f^{'}(0)=\left.\frac{df(x)}{dx}\right|_{x=0}$.

By \eqref{23}, we get
\begin{equation}\label{24}
q\beta_{n,q}(x+1|\lambda)-\beta_{n,q}(x|\lambda)=(q-1)\lambda^n\left(\frac{[x]_q}{\lambda}\right)_n+\sum_{l=1} ^{n-1}S_1(n,l)\lambda^{n-l}l[x]_q ^{l-1}q^x,
\end{equation}
where $n\in{\mathbb{N}}$.

Therefore, by \eqref{24}, we obtain the following theorem.
\begin{theorem}\label{thm9}
For $n \geq 0$, we have
\begin{equation*}
q\beta_{n,q}(x+1|\lambda)-\beta_{n,q}(x|\lambda)=(q-1)\lambda^n\left(\frac{[x]_q}{\lambda}\right)_n+\sum_{l=1} ^{n-1}S_1(n,l)\lambda^{n-l}l[x]_q ^{l-1}q^x.
\end{equation*}
\end{theorem}


\begin{thebibliography}{10}

\bibitem {01} A. Bayad, J.  Chikhi, {\it Apostol-Euler polynomials and asymptotics for negative binomial reciprocals}, Adv. Stud. Contemp. Math., ${\mathbf{24}}$ (2014), no. 1, 33-37.

\bibitem {02} J. Choi, T. Kim, Y. H. Kim, {\it A note on the extended $q$-Bernoulli numbers and polynomials}, Adv. Stud. Contemp. Math., ${\mathbf{21}}$ (2011), no. 4, 351-354.

\bibitem {03} L. Carlitz, {\it Degenerate Stirling, Bernoulli and Eulerian numbers}, Utilitas Math., ${\mathbf{15}}$ (1979), 51-88.

\bibitem {04} L. Carlitz, {\it A degenerate Staudt-Clausen theorem}, Arch. Math. (Basel), ${\mathbf{7}}$ (1956), 28-33.

\bibitem {05} D. Kang, S. J. Lee, J.-W. Park, S.-H. Rim, {\it On the twisted weak weight $q$-Bernoulli polynomials and numbers},   Proc. Jangjeon Math. Soc.,  ${\mathbf{16}}$ (2013), no. 2, 195-201.


\bibitem {06} T. Kim, {\it $q$-Volkenborn integration}, Russ. J. Math. Phys., ${\mathbf{9}}$ (2002), no. 3, 288-299.


\bibitem {07} J. W. Park, {\it New approach to $q$-Bernoulli polynomials with weight or weak weight}, Adv. Stud. Contemp. Math., ${\mathbf{24}}$ (2014), no. 1, 39-44.

\bibitem {08} J.-J. Seo, S.-H. Rim, S.-H. Lee, D. V. Dolgy, T. Kim, {\it $q$-Bernoulli numbers and polynomials related to $p$-adic invariant integral on ${\mathbb{Z}}_p$}, Proc. Jangjeon Math. Soc.,  ${\mathbf{16}}$ (2013), no. 3, 321-326.

\end{thebibliography}
\end{document}